# Emergence of Mathematics in Ancient India

## A Reassessment

Jaidev Dasgupta

Independent Researcher, jaidevd101@gmail.com

**Abstract**

This work explores a possible course of evolution of mathematics in ancient times in India when there was no script, no place-value system, and no zero. Reviewing examples of time-reckoning, large numbers, sacrificial altar-making, and astronomy, it investigates the role of concrete objects, natural events, rituals and names in context-dependent arithmetic, revealing its limited scope confined to counting, addition and subtraction. Higher operations, namely, multiplication, division and fractional calculations had to wait until the advent of symbolic numerals and procedures for computation. It is argued that the impression of these higher operations in a period usually known as the Vedic times is caused by inadvertent interpolation of present knowledge of mathematics in modern readings of the ancient texts.

**Introduction**

Various theories have been proposed about the antiquity, the roots, and the advancement of Indian civilization (Thapar, 2004; Witzel, 2001; Bryant, 2001). Among others, the knowledge of mathematics has also been considered an important mark of the civilization and the level to which it was advanced. Mathematics being a study of patterns, order and regularity, underlies the knowledge of science and technology in a culture. Considering Harappan civilization at its height in the middle of the third millennia BCE, there is ample evidence supporting their knowledge of mathematics, measurement and technology in the rectangular grid roads, public baths, buildings, granaries and residences. But after its decline in 1900-1700 BCE, the subsequent period, which is usually called the Vedic period, (1500 BCE-500 BCE) does not have evidence in terms of archaeological artifacts, buildings, or rock inscriptions to demonstrate such knowledge. However, the texts from that period have extensive references to mathematics - especially in astronomy, calendrics for time-reckoning, and religious, sacrificial altar-making. The *Sulva Sutras* for alter-making and *Vedanga Jyotish*, the astronomical text for calendar making, are prime examples of such knowledge in antiquity. Seidenberg's remarks about rituals as the source of knowledge of mathematics in ancient cultures (1981) has further spurred a wave of discovering arithmetic and astronomical knowledge in ancient Indian texts to establish a greater degree of mathematical knowledge in the remote past.

This search for mathematics in antiquity, however, has overlooked the question of the nature of mathematics and its scope and feasibility in the particular context of Vedic culture when there was no script (hence no numeral system), no place-value system, and no zero. The question that deserves greater attention is the emergence and evolution of mathematics in that period with those constraints. Instead, the current knowledge of mathematics often gets projected back into the several millennia old texts that are believed to have come down through oral means. These attempts to explain ancient mathematics using modern tools, symbols, and systems of mathematics, often cause puzzlement over how such calculations were possible in a period without these tools and techniques, and to what extent the claims of accuracy and achievements are supportable. The present work attempts to explore some of these questions from a practical and mathematical points of view, and tries to gauge the kind and level of arithmetic that existed in the ancient times.

**Hypothesis**

It has been suggested that a culture may continue to survive in the absence of a writing system (Kenoyer, 2020). Which may well be the case, but this absence may leave its stamp on the nature and capabilities of the culture. Strong emphasis on memorization and orality of the Vedic tradition can certainly be related to the lack of a script. While this lacuna led to the development of poetic meters or *chandas* and composition of hymns following these metric rules for recording, communicating, and transmitting knowledge, it impacted other aspects of the culture. For instance, mathematics that heavily depends on external symbolic representation of numbers and generation of procedures for manipulating them, remained constrained in this period due to the absence of a writing system. How numbers were visualized, represented and generated, how time was tracked and astronomical calculations performed for calendar-making, or how calculations tied to altar-making were conceived are some of the questions that deserve closer scrutiny to reach a better understanding of mathematical capability in the Vedic culture. The following discussion shows that arithmetic of that period was limited in its scope and capability as it was context-dependent, tracking and tallying concrete objects, natural events, rituals and names.

**Discussion**

*Time-Reckoning in Rituals and Altars*

Across different cultures one finds a common tendency of tracking time. The Vedic texts also display the same spirit. Starting from the day as a natural unit of time between two consecutive sunrises or sunsets to figuring out the number of days in a month and in the two fortnights (*parvas*) of a month, the number of seasons, the twelve months in a year with 360 days, and resolution of the length of a year between 365 and 366 days - all reveal pattern recognition, tracking of a natural order, and the associated mathematical activity of counting and record keeping of elapsed days by following annual rituals.

A year-long ritual, the *gvamāyana sattra*, was conducted in which oblations were offered every morning and evening, every new and full moon (i.e., *parvas*), and at the commencement of every season and *ayana*, when the direction of solar movement changes from north to south and vice versa at the two solstices (Satapatha Brahman SB 1.6.3.35-36; Tilak, 1893). These periodic sacrifices may be considered to have served the same purpose as the tally marks on a stick or a piece of bone in prehistoric times used for time-tracking by following lunar phases (Marshack, 1991; Livio, 2002).

The concept of *sadaha* was yet another tool for tracking time by tallying or counting, in which a set of six consecutive days were identified by specific rituals (Taittiriya Sanhita TS 7.2.1; SB pt. III, p. xxi; SB pt. V, p. 148; Dikshit, 1969). Sadaha divides a thirty-day month into five equal parts. Every time the full moon returned on the $29^{th}$ day instead of $30^{th}$, the last day of the sadaha was dropped to keep in sync with the lunar phase.

Thus, the annual ritual performed the function of a calendar by tracking solar and lunar movements. It entailed counting number of days as one continued through fortnights, months, seasons and a year, one day at a time, dropping a day here, adding another there, using no more than simple addition and subtraction to reconcile the lunar and solar years. A sunrise, a sunset, a lunar phase, including the full- and the new-moon served as tokens or counters meant for external representations of numbers.

Another means of tracking time is found in the extensive naming scheme used for years (five names of years in a five-year *yuga*), the twelve months, the dark (*krishna*) and bright (*shukla*) fortnights or *parvas* of a month, the names for 24 parvas in a year, the two parts of each day in the light and dark halves in a month adding up to sixty distinct names of days (*ahas*) and nights (*rātris*) in a month, the 30 divisions of a day called *muhurtas*, each with its own name in the two fortnights or *paksha* in a month, and each muhurta further subdivided into fifteen *prati-muhurtas* with their own specific names (Dikshit, 1969). These distinct names represented specific segments of time and provided 810,000 (5 years x 12 months x 2 parvas x 30 day/night(*ahorātra*) x 15 muhurtas x 15 prati-muhurtas) unique sequences of names pointing to as many specific time segments of prati-muhurtas (3.2-minute long) in a five-year yuga.

Whether such long sequences of names were ever used is not known, though, in principle, they had the potential to be used when numerical representation of time such as 12/18/2021, 15:51 hr. – for December eighteen, year two-thousand-and-twenty-one, at three-hours fifty-one minutes in the afternoon – was not invented (due to lack of script). Ancient texts (and inscriptions in later periods) indicate use of shorter combinations of fewer elements – e.g., month, *paksha*, *tithi*, muhurta, and *nakshatra*. *Tithi* is a lunar day (a variable unit of time between 20 and 26.8 hrs. with average of 23.62 hrs.) based on the time taken by the moon between two consecutive moon rises or sets (Saha and Lahiri, 1955, p. 221). Each half of a month has fifteen tithis with ordinal names: *Prathama* (first), *dvitiya* (second), *tritiya* (third), *Chaturthi* (fourth) and so on until the fifteenth day, called either *purnamasi* (full moon) or *amavasya* (new moon), depending on the bright or the dark half of a month, respectively. Even *nakshatras* or star names were used for this purpose. In short, again, in the absence of numeral symbols, names – which represent or identify a thing or an entity – were used to point to specific times. The practice of using some of these names is still prevalent in India, especially in astrological and religious contexts.

Using names for tracking time does not mean the ancient Indians were unfamiliar with numbers. Verses (1.164.10-15, 48) in the Rig Veda (RV), for example, mention numbers three, five, six, seven, ten, twelve, three-hundred-and-sixty, and seven-hundred-and-twenty. One of these verses, no. 11, speaks of a year with 720 days (*ahas*) and nights (*rātris*) joined in pairs, together suggesting a year of 360 days. The two verses RV 10.62.7-8 speak of a thousand (*sahasra*) cows and a hundred (*śata*) horses. The Taittiriya Samhita sections 7.2.11-20 and the use of poetic meters or *chandas* also offer ample evidence of the knowledge of numbers.

However, there is no information available on number symbols and their use in mathematical operations other than a slim clue to the use of a numeral in a verse 10.62.7 of RV, which is open to interpretation: According to one translation of the verse, number eight is marked on the ears of a thousand cows (Griffith, 1896), whereas another translation simply says "cut-branded ears" (Jamison and Brereton, 2014). Aside from this single, weak and isolated reference, the extant ancient literature on numerals is silent on symbolic representation of numbers, which makes perfect sense considering there was no writing system at that stage that used representative symbols. But there is significant evidence of counting, creating number words or number names and some arithmetical operations in that period.

According to Murthy (2005), the Vedic bards were fascinated by early mathematics and associated spiritual or religious symbolism with small numbers. They were comfortable with simple arithmetic, e.g., addition. In the case of larger numbers, such as three-hundred-and-sixty (number of days in a year), he goes as far as to state that "I do not wish to write it as '12 x 30' or '24 x 15' or '3 x 3 x 4 x 10', as I am not sure that the Vedic people are comfortable with multiplication that results in a three digit number." Instead, he proposes to write the number as summation of smaller numbers.

But the Shatapatha Brahman (SB), a text of around middle of the first millennium BCE (Witzel, 2001), seems to display the knowledge of numbers and multiplication. Sections SB 10.4.2.4-17 are said to present numbers 2, 3, 4, 5, 6, 8, 9, 10, 12, 15, 16, 18, 20 and 24 as the divisors of 720 - the 360 pairs of day and night mentioned above (Kak, 1993). Number 1, the fifteenth divisor, is implied in the section 2 which says that "…in this Prajapati, the year, there are seven hundred and twenty days and nights, his lights, (being) those bricks…." The Brahman uses Prajapati as equivalent to a year.

Caution must be exercised here against interpolating modern ways of making 720 by multiplying its factors. The proposed set of 15 numbers as the divisors or factors of 720 in SB, and hence of the evidence of multiplication, demand an explanation as to how in the absence of number symbols and operations between them this could have been possible. Instead, it would be more appropriate to consider that the text actually enumerates possible ways of distributing or splitting 720 bricks into smaller collections of bricks. Each verse has the form: He (*Prajapati*) made himself four bodies of a hundred and eighty bricks each; … He made himself six bodies of a hundred and twenty bricks each; …He made himself ten bodies of seventy-two bricks each; and finally, …He made himself twenty-four bodies of thirty bricks each (SB 10.4.2.4-17). The repeating phrase "made himself" provides the clue to creating collections of concrete objects, such as bricks, and counting and adding them; which was the practical way of keeping account of things prior to the invention of written symbols of numbers and abstract procedures for addition and multiplication, neither of which find any mention in the Vedic texts.

To have a better appreciation of how bricks were used for representation, let us look at a few examples. In SB 10.4.2.30 an enclosing-stone of an altar is equated to a night (*rātri*) and its 15 muhurtas, and a *yajushmati* brick is equivalent to a day (*aha*) and its 15 muhurtas; thus, making three-hundred-and-sixty each of enclosing-stones and *yajushmati* bricks in an altar (SB 10.5.4.5). Elsewhere, in the same Brahman, a *yajushmati* brick also refers to half-moons (fortnights), months, and seasons; and *lokamprina* bricks to the muhurtas in a year (SB 10.4.3.12). Altar bricks refer to the nakshatras or asterisms as well (SB 10.5.4.5), that are markers of time in their own turn. Thus, the abstract notions of time points and units are tied to material objects, which indicates an earlier stage of mathematics prior to representing the abstract by written symbols.

It must be mentioned, here, that while a stone or a brick is used as a unit pointing to a bundle of time units, e.g., muhurtas or days, there is a lack of standardization (more instances of it later). The *yajushmati* brick is used to refer to multiple entities, and both a stone and a brick refer to the muhurtas. At the same time, the qualitative difference between day and night is preserved by using a brick for the day and a stone for the night. Thus, the meanings of these concrete objects are fixed in the context of a fire-altar which is to be constructed following an order to represent a year. Use of unlimited number of stones and bricks or making altars of arbitrary size are forbidden (SB 10.4.3.5-6).

Thus, the Shatapatha Brahman is pointing to the use of more ancient procedures of representation and computation than what is usually believed. To contrast it from arithmetic that uses abstract numerals, it may be termed proto-arithmetic which is embedded in the context of a social activity and tied to concrete objects. This is what Seidenberg meant with his remark that the roots (or inception) of mathematics are to be found in social rituals. With altar-making and fire sacrifices as characteristic social processes of the Vedic culture, their use of stones and bricks in altars to refer to mathematical entities was quite inventive.

The text itself furnishes ample evidence of how natural patterns of annual changes in time and seasons were tracked and computed, using names, rituals, and concrete objects as representations or symbols in the absence of written numerals. This was a significantly advanced stage of mathematics with a larger scope and capability as compared to the mathematics of prehistoric times confined to counting one, two, and three… and maybe up to ten on fingers.

Nevertheless, concrete representations have their own limitations. Bricks and stones allow physical operations of collecting, combining, removing, sorting, separating, distributing, apportioning, and counting, that are quite different from the mathematical operations performed on numerals. While the former, as remarked above, is fixed and embedded in its context with a specific purpose and meaning, the latter is independent of it, permitting freedom for growth by manipulating written symbols. Additionally, the altar-based mathematics was not standardized either. Hence, the nature of mathematics was quite different in the Vedic period from what it is now since the second half of the first millennium. With the advent of script and written numerals the former evolved into the latter.

*Large Numbers*

Besides small numbers encountered in early astronomy, sacrifices and calendrics, often, the evidence of mathematical prowess of the ancients is presented in their ability to generate large numbers. White (*Shukla*) Yajurveda, for instance, mentions thirteen number names:

*eka* (1), *daśa* (10), *śata* (100), *sahasra* (1,000), *ayuta* (10,000), *niyuta* (100,000), *prayuta* (1,000,000), *arbuda* (10,000,000), *nyarbuda* (100,000,000), *samudra* (1,000,000,000), *madhya* (10,000,000,000), *anta* (100,000,000,000), *parârdha* (1,000,000,000,000) (Datta and Singh, 2004).

These number names are lifted from YV17.2 without their decimal representations as presented above in the parentheses:

" [imā́](#) [me](#) [agna](#) [íṣṭakā](#) [dʰenávaḥ](#) [santv](#) [ékā](#) [ca](#) [dáśa](#) [ca](#) [dáśa](#) [ca](#) [śatám](#) [ca](#) [śatám](#) [ca](#) [sahásram](#) [ca](#) [sah ásram](#) [cāyútaṃ](#) [cāyútaṃ](#) [ca](#) [niyútaṃ](#) [ca](#) [niyútaṃ](#) [ca](#) [prayútaṃ](#) [cā́rbudaṃ](#) [ca](#) [nyàrbudaṁ](#) [samudráś](#) [ca](#) [mádʰyam](#) [cā́ntaś](#) [ca](#) [parārdʰáś](#) [caitā́](#) [me](#) [agna](#) [íṣṭakā](#) [dʰenávaḥ](#) [santv](#) [amútrāmúṣmiml](#) [loké](#) "
(TITUS: YVW, 2012)

And here is the translation of this verse by Griffith (1899):

"O Agni (fire), may these bricks be mine own milch kine: one, and ten, and ten tens, a hundred, and ten hundreds, a thousand,

and ten thousand a myriad, and a hundred thousand,

and a million, and a hundred millions, and an ocean middle and end, and a hundred thousand millions, and a billion.

May these bricks be mine own milch-kine in yonder world

and in this world."

A few other Vedic texts, namely, the Taittiriya Samhita (TS 4.4.11), Kathak Samhita (KS 17.10) and the Maitriyani Samhita (MS 2.8.14) also repeat the same verse either with the same sequence of number names or slightly different ones.

Two important questions worth considering here are: How were these numbers generated? And what did they signify or mean?

As the first four number names in the Yajurvedic series – eka, daśá, śáta and sahasra – that are multiples of ten and are still in use with meaning one, ten, a hundred and a thousand respectively, it may be assumed that the whole series of 13 numbers is based on the decimal system. Since the decimal system is believed to have arisen from counting on ten fingers in the Upper Paleolithic/Mesolithic period, it is reasonable to assume that in the middle of 2$^{nd}$ millennium BCE the system was already in use in India. Egyptians were also using 10-base system around 3000 BCE.

This leads to the question of how these numbers were created in the absence of symbolic representation when concrete objects were in use to relate to numbers. Interestingly, similar to the Shatapatha Brahmana, the Yajurveda too uses concrete objects to refer to the numbers: "*May these bricks be mine own milch-kine* (cows)" (see above, Tr. by Griffith). The other texts, named above, also use the same reference. Starting from smaller number of bricks and cows – one, ten, a hundred, and a thousand – they proceed to conceive larger and larger numbers of bricks and cows. These texts are consistent with the contemporary cultural practice of using concrete objects for representing numbers.

To further deliberate on the generation of long numbers and their significance, we need to review a few other numbers from the period several centuries after the end of the Vedic age (~500 BCE). Taking two examples from the Jain and Buddhist texts of around the turn of the era. *Lalitavistāra*, a Buddhist text, (Ch. 12) gives a series of 24 number names that proceeds on a centesimal scale starting from one-hundred times *koti* (or ten million) and stepwise reaches the number $10^{53}$, named *tallakṣana*. Similarly, a Jain text *Anuyogadvara-sutra* mentions a 29 digits long number $2^{96}$ obtained by multiplying the sixth square of two by its fifth square, i.e., (18,446,744,073,709,551,616) x (4,294,967,296) (AYD, pp. 384-385).

It is to be noted that these incredibly long numbers do not suddenly appear in one's mind, they are generated starting from smaller numbers either by multiplying them by 10 or a hundred or by squaring them repeatedly and, often, at each step a name is generated to keep track of the numbers being produced. Recall the use of names as markers for referring to time segments (above). Use of such mechanical procedures for producing large numbers could simply be to display mental and mathematical prowess of the person producing them, such as bodhisattva in Lalitavistāra. But the ability to generate large numbers mechanically and understanding their significance or meaning are entirely different things.

The Jain text claims that $2^{96}$ is the population of human beings in the world. And Lalitavistāra says that tallakṣana ($10^{53}$) is the basic unit for calculating the size of Mount Meru. Furthermore, according to bodhisattva, the next larger number beyond tallakṣana is *dhvajāgravatī* ($10^{55}$), using which as a unit, it is possible to calculate the grains of sand in the river Ganges. Mount Meru is a mythical mountain measuring the height of which is as meaningless as counting the sand grains in the river Ganges. Similarly, an assumption of the human population of $2^{96}$ (~ $10^{29}$) on earth in the Anuyogdvara sutra, composed around 100 BCE (Datta and Singh, 2004), is equally fantastic and unreal; one wonders how the author reached that conclusion. Therefore, although the reference to sand grains, Mount Meru or to the world population provide a context for the mind to grasp these numbers, the best that can be claimed – in the light of their unreality and sheer extravagant size – is that by applying a formula or algorithm repetitively any large number, howsoever long, could be generated and named; and that this art was known to people in India since a long time.

Perhaps such exercises were indulged in only to establish that starting from a small number one can reach, in a limited number of steps, such immense numbers that are simply incomprehensible and beyond human imagination and experience, conveying a sense of infinitude. Generation of long numbers starting from *eka* or 1 to *parârdha* (1,000,000,000,000) in YV17.2 was perhaps one of the earliest such attempts. As Vedic philosophy is well-known to be imbued with the thoughts of infinite and infinity, it may be conjectured that the sequence of long numbers in a religious, ritualistic text like Yajurveda was perhaps meant to guide one's mind out of the practical world, stepwise, into the metaphysical world.

The last comment does not mean that such sequences were confined solely to the Vedic religious texts. As already discussed, the Jain and the Buddhist traditions too had similar large numbers. Fascination for long numbers persisted in mathematical works well after the first millennium of the Common Era. A notable point about these number names is that often they don't seem to represent a fixed value. While ten-, a hundred- and a thousand-millions are *arbuda*, *nyarbuda* and *samudra* in Yajurveda17.2, in the 5$^{th}$ century Aryabhatia (2.2), the same values are named – *koti*, *arbuda* and *vraṇda*, respectively (Shukla and Sarma, 1976). Similarly, while an *ayuta* is ten-thousand in the Yajurveda, in Lalitavistāra it is a billion. Even after Indian mathematics had already made major advances in the first millennium of this era, the lists of numbers quoted by Sridhara (8$^{th}$ century), Mahavira (9$^{th}$ century), Bhaskar II (12$^{th}$ century), and Narayana (14$^{th}$ century) are often different in their names and values (Datta and Singh, 2004, pp. 12-13). Al-Biruni in his book on India also draws attention to such variations in names of numbers (Sachau, 1910, pp. 175-177). Thus, while the tradition of generating and naming large numbers started from the Yajurvedic time, it lacked a standard set of names with fixed magnitudes other than the first four terms – eka, daśa, śata and Sahasra – even after more than two millennia. It is therefore possible that the names of large numbers in the Vedic period were not supposed to have any significance other than being the tokens for memorization of mechanically generated numbers.

On the other hand, the constancy of the use of the first four names – eka (1), daśa (10), śata (100) and Sahasra (1,000) – since the time of the Rig Veda indicates that people in that period were dealing with these numbers frequently in daily commerce, astronomy, and time-reckoning, and hence were familiar with their names and values, making them standard terms. It is in the *Arthashastra* (Shamasastry, 1915; Rangarajan, 1992), compiled between 300 BCE to 300 CE, that one comes across larger numbers in practical life when the salaries of state employees were in the order of tens of thousands of *panas*, and fines for certain offences were as high as a few thousand *panas*. That's when people could relate to numbers higher than sahasra or a thousand. Prior to that, in the agropastoral economy of the Vedic times, sahasra was most likely the largest number people came across. Hence the doubt expressed by Murthy: if the Vedic poets had "…an absolute idea of how big the number '60099' is…" (Murty, 2005), which is the number of men killed by the mythic god Indra in a battle (RV1.53.9).

Thus, while exploring the significance of large numbers such as a million, a billion or higher in the ancient period, one must ask about their presence and use in social life. According to the Material Engagement Theory in anthropology, counts extend as far as the number of resources and material possession in a culture (Barras, 2021). For instance, ancient Egyptians had counts up to a million; they counted 400,000 oxen, 1,422,000 goats, and 120,000 prisoners as gains from the victory of the predynastic king Narmer over lower Egypt (Clagett, 1999; p. 2). Likewise, cows were the mark of wealth in the Vedic pastoral economy. Again, it is the concrete which helps grasp the meaning of a number. Following the period of the Arthashastra, the next phase when numbers much larger than ten or hundred thousand assumed significance was in the second half of the first millennium in performing astronomical calculations.

*Astronomy*

Vedanga Jyotish (VJ), which is considered the astronomical treatise of the Vedic times, is yet another text providing a window to the knowledge of mathematics of that period. It was used for calendar making. While usually the discussion on VJ revolves around its astronomy and the time of its composition – whether 1400 or 1200 BCE or 400 BCE – rarely ever is the feasibility of mathematics in the text in the proposed age a matter of investigation. It is assumed that the required knowledge of mathematics for astronomical computations was already available.

Though a detailed presentation of VJ is out of scope here, to set the context for the discussion below, here is a brief picture of it: Vedanga Jyotish text uses a model of a 5-year yuga cycle in which the cycle begins with the conjunction of the sun, the moon and the Dhanishthā nakshatra or asterism at winter solstice, which repeats every five years with the beginning of a new cycle. The sky in VJ is fixed, there was no knowledge of precession at that time. The solar year is considered 366 civil days long, yielding 1830 civil days in a yuga. In this period, there are 62 lunar synodic months (29.5 days long) and 67 lunar sidereal months (27.3 days long). For calendar making purposes, an important question is to determine the nakshatras in which the full moon and the new moon will occur during the quinquennial yuga.

Since in a yuga there are 67 sidereal months during each of which the moon traverses through 27 nakshatras, in a five-year period it traverses through 1809 (67x27) asterisms. At the same time, the moon also goes through 124 fortnights (2 fortnights in each of 62 synodic months). Hence, the moon covers 1809/124 (= 14.58871) asterisms per fortnight. Obviously, there is fractional calculation involved here. To write the same number of nakshatras traversed per fortnight in a non-decimal form it is 14 73/124. As the moon keeps moving from one fortnight to another, the same number with its fractional part gets added repeatedly.

An important question to ask here is: In the absence of numerals, the place-value system and zero how much of the above calculation possible in the Vedic period? These calculations also involve knowledge of dealing with fractions. Several other calculations in VJ depend on fractional arithmetic. For instance, the rule of three involves knowledge of fractions (ratios), their equality and the derivation of an unknown quantity, a, from three known quantities b, c, and d (a/b = c/d; therefore, a = b.c/d). This derivation is an algebraic procedure. Was the knowledge of mathematics advanced enough in that period to support such procedures and calculations? For modern readers, these mathematical operations and the systems enabling them are so obvious that they rarely ever stop to question their feasibility in the ancient texts. Anachronistically, they project their knowledge of mathematics into antiquity.

It is well known that the decimal place-value notation and a symbol for zero appeared in Indian mathematics between $5^{th}$ to $7^{th}$ century CE (Katz, 2009; Joseph, 2016). Inscriptions, epigraphs and grants from the $3^{rd}$ century BCE – when the Brahmi script first appeared in the Indian subcontinent

– till the first few centuries of CE show that single symbols stood for denoting the numbers 1, 10, 100, and 1000 (see more below). Given these historical facts, the mathematics in the extant version of VJ seems inconsistent with the age it is usually portrayed in. For it to be a text of the Vedic time, it is important to establish the feasibility of its mathematics, especially of fractional calculations, in that period.

The knowledge of fractions in ancient India is usually believed to exist since the time of the Rig Veda as the *Purusha* hymn mentions $1/4^{th}$ and $3/4^{th}$ of Purusha or the Original Being (RV 10.90.3-4). In fact, the hymn simply says that – one foot (*pāda*) of the Being is in this world and three feet (*tri-pāda*) above in the heaven (*pādo.asyaviśvā bhūtāni tripādasyāmṛtaṃ divi*; also see Trivedi, 1954) – which is in agreement with the age when numbers were still dealt with the help of concrete objects such as bricks, stones, fingers (*angula*) and feet (*pāda*). Although at some later point in time, pāda also acquired the meaning of one-quarter, there is no way of establishing that RV means fractions. This could be due to the present familiarity with fractions which leads to reading one-fourth and three-fourths in RV - a clear example of how present knowledge gets interpolated in the ancient texts.

Even if it is assumed that pāda meant one-of-four or one-fourth in RV, there is a big difference between having a simple knowledge of parts or fractions of a thing and performing mathematical operations on these fractions. The problem surfaces when one has to add 1/3 and 3/4, or divide 2/3 with 4/5. Procedurally, operations on fractions aren't as easy or intuitive as it seems, they require knowledge of proper methods.

Modern research on evolution, teaching, and learning of mathematics shows that though mathematics is a language, it is unlike a natural language that people pick up easily in childhood; arithmetic has to be learnt which includes learning of language, numbers, symbols and operations. Fractions and mathematical operations on them are even harder to understand than dealing with whole numbers (Dehaene, 1997; Devlin, 2000; Caswell, 2007; Heath and Starr, 2022). People can readily relate to and develop a feeling for integral numbers rather than fractions. Being a part of something, fractions pose the problem of how to deal with them unless one knows the procedure of using the least common denominator, which is not an easy concept and requires prior knowledge of multiplication and factorization. Moreover, without symbolic numerals, performing such mathematical operations are incomprehensible. The silence of the Vedic texts on methods for arithmetical operations reflects their inability to conduct them in an age when numbers were conceived in terms of concrete objects.

In this regard, the *Arthashastra* – being chronologically situated between the Vedic period (1,500-500 BCE) and the classical period of Indian mathematics (400-1600 CE) – might shed some light on the status of mathematics in that period, which could, in turn, help address the question under consideration. It is interesting to note that although the book is an encyclopedic compilation of a wide range of knowledge, it neither mentions Vedanga Jyotish nor says anything about the lunar mansions or nakshatras so integrally tied to the VJ astronomy; though it describes the five-year yuga, the equinoxes and the solstices, six seasons, twelve months, the synodic and sidereal lunar months, the two fortnights, the intercalary months and when they are added. But surprisingly, the astronomical conjunction of the sun, the moon, and the Dhanishthā asterism at winter solstice with which the five-year yuga begins in VJ is missing. This suggests that at that time the sidereal calculations of VJ were unknown.

Furthermore, though the text lists units of weights, lengths, measures of capacity and divisions of time, it does not say anything about mathematics of its time. While it may be argued that a book on statecraft is not expected to cover mathematics, a total silence on this topic is also intriguing especially since it is a book of economics as well which entails calculations. However, unlike other texts prior to it, the *Arthashastra* extensively mentions fractions – predominantly of one – in the description of taxes, tolls, fee, surcharges and service charges, indicating their use in practical life. At the same time, different sets of weights and balances used for different purposes and other measures were so well defined in whole numbers that a daily use of fractional computation is doubtful. Even in the case of fractional taxes in kind there was no need for fractional calculation. For example, $1/6^{th}$ tax on agricultural production could be paid simply by separating the sixth measure (for tax) after every five measures (by weight of any unit). Nonetheless, their familiarity with computation using unit fractions cannot be ruled out. Hence, the *Arthashastra* shows a distinct step forward from knowing unit fractions to using them in calculations.

It must be mentioned here that the *Sulva Sutras* – the texts for constructing sacrificial altars for *yajna* – that are dated prior to as well as contemporary to the *Arthashastra,* also appear to have the knowledge of fractions (Thibaut, 1875; Dani, 2010). However, the way numbers are tied to bricks and stones in altars (see above), the geometry used in designing and constructing of these altars is also practice bound and concrete – using ropes (*rajju*) and pegs and bamboos. It is one thing to extend or reduce the length of a rope by a part of it to draw a circle or a square or other geometrical figures of a certain area, it is yet another to do fractional calculations using symbolic numerals. While the former is concrete, the latter is abstract. Since the Brahmi script appears in the $3^{rd}$ century BCE, it is hard to believe the existence of any symbolic manipulation of numbers and fractions prior to that time.

Even after that, it remains to be explored how these mathematical operations were possible with symbols in Brahmi script. As a new script, it took some time to develop a full set of notations. The data collected from coins, inscriptions and grants from the first and second century CE onwards show the use of individual symbols for numbers 1 through 10, multiples of 10, of 100, and even of 1000. Brahmi numerals neither had a symbol for zero nor was it a system with positional value. These numerals were in use right up to the $6^{th}$ or $7^{th}$ century CE (see Tables III-XI in Datta and Singh, 2004).

Which takes us back to the question of feasibility of sidereal and other fractional calculations in Vedanga Jyotish. When were the required methods for computations developed? As remarked above, fractions in VJ are composite and more complex to deal with than the unit fractions in the *Arthashastra*. According to Datta and Singh (2004), the earliest record of reduction of fractions to lowest terms (entailing the knowledge of multiplication and factorization) is in a Jain text *Tattvarthadhigham Sutra*, which was believably written sometime between the $2^{nd}$ and the $5^{th}$ century CE. The next authors to write on operations on fractions were Brahmagupta in the $7^{th}$ century and then Sridhar in the $8^{th}$.

Taken together, a few possibilities emerge: one, that the five-year yuga cycle and the associated synodic lunar monthly calculations were already known prior to VJ; two, the *Arthashastra* was composed before the extant VJ, otherwise, there is no reason for the former not to mention the latter and miss important information on moon's sidereal locations, especially when the court employed astrologers whose job was to compute for calendrical purposes; and three, the extant VJ was composed later by combining the pre-existing synodic calculations of the moon with its positions in the nakshatras when appropriate mathematics became available. The knowledge of the conjunction of the sun, the moon, and the Dhanishthā nakshatra at the winter solstice transmitted orally over centuries was taken as a reference starting point for calendrical computations. With a fixed sky (as VJ does not mention precession of the equinoxes), these computations could be performed anytime later.

**Conclusion**

While looking back at the numbers and calculations in the past, modern scholarship often reflexively projects a mathematical system with which they are familiar. Literature is filled with examples of understanding ancient Indian math in modern terms, and drawing conclusions on that basis paying little attention to the fact that in that remote age there were no symbolic representation of numbers, a place-value system, or zero. If we ask the question of how the ancients performed mathematics without these concepts and systems, the limitations under which they operated would become obvious. Using examples of practices in the Vedic period, this article has attempted to reveal such constraints and suggested possible ways in which some of these calculations were performed, for instance, using repeated addition instead of multiplication. While multiplication and division are easier and faster than repeated addition or subtraction from the modern perspective, ancients had no recourse but to go through the routine of addition and subtraction by counting concrete objects or tallying natural events or tracking rituals. Even the size of numbers encountered in that period were most likely no larger than a thousand, other than the mechanically generated large numbers with doubtful understanding of their significance. Dealing with fractions itself was another roadblock. Understanding fractions in terms of parts of something concrete, say, a rope or a brick or agricultural product, is different from writing and manipulating fractions in symbolic form. Only after the advent of a script and other systems, including mathematical operations to deal with fractions, the calculations in the extant Vedanga Jyotish seem to have become possible. Overall, from this line of reasoning emerges a path of slow evolution of arithmetic and geometry, which helps with a better appreciation of the history of mathematics in ancient India. It tells more about the ingenuity of the ancients, how they accomplished their tasks and goals with limited mathematical resources; and once the required systems were in place during the first few centuries of the Common Era, Indian mathematics leaped forward rapidly from the 5[th] century onward.